\documentclass[11pt]{article}
\usepackage{amsmath,amsthm,amsfonts,amssymb,latexsym}
\usepackage{tikz}\usetikzlibrary{intersections}
\usepackage{graphics}
\newcommand{\pict}{%
\begin{center}
\begin{tikzpicture}[scale=1.0]
    \draw (0,0) node [anchor=south east] {\(z_t\)} -- (3,0) node [anchor=south] {\(v\)} -- (-1,-3) node [anchor=north] {\(w\)} -- cycle;
    \draw (0,0) -- (-2,0) node [anchor=east] {\(u\)} -- (1,-1.5) node [anchor=south east]{};
    \draw (1,-1.5) node [anchor=north west] {\(v_s\)};
    \draw (-1,-0.8) node [anchor=north west] {\(z_{t_s}\)};\end{tikzpicture}
\end{center}
}%

\newcommand{\R}{{\mathbb R}}

\newcommand{\la}{\langle}
\newcommand{\ra}{\rangle}

\newtheorem{theorem}{Theorem\hspace{-0.4em}}

\begin{document}
\thispagestyle{empty}
\begin{center}
{\Large\bf Linear Perturbations of Quasiconvex Functions \smallskip\\
and Convexity}
\end{center}
\begin{center}
{\large  Khanh, Pham Duy\footnote{Department of Mathematics, University of Pedagogy of Ho Chi Minh City, 280 An Duong Vuong,
Ho Chi Minh, Vietnam, \textbf{\tt pdkhanh182@gmail.com}} and
 Lassonde, Marc}\footnote{Universit\'e des Antilles et de la Guyane,
 Campus de Fouillole - BP 250,
  97157 Pointe \`a Pitre, France, \textbf{\tt marc.lassonde@univ-ag.fr}}
\end{center}

Let $E$ be a real vector space
and let $C\subset E$ be a nonempty convex subset.
We recall that a function $f : C\to\R$ is said to be
{\it convex} if for all $u,v\in C$ and $t\in {]}0,1{[}$,
$$
f(v + t(u-v))\le f(v)+ t(f(u)-f(v)),~~
$$
and {\it quasiconvex} if for all $u,v\in C$ and $t\in {]}0,1{[}$,
$$
f(v + t(u-v))\le\max\{f(u), f(v)\}.
$$

It is well known since the pioneering work \cite{Cro77},
and easy to prove (see e.g. \cite[Proposition~2.1]{ACL94}),
that a function $f$ is convex provided all its linear perturbations
$f+u^*$, $u^*\in E^*$, are quasiconvex.
The purpose of this note is to show that
if a function $f : C\to\R$ satisfies a mild stability property
at 'flat' points of the (relative) boundary of $C$,
the convexity of $f$ is guaranteed as soon as
for some $c^*\in E^*$ not constant on $C$, the function
$f+\lambda c^*$ is quasiconvex for all $\lambda\in\R$.
\medbreak
We say that a function $f:C\rightarrow \R$ is
\textit{radially lower stable at $z\in C$},
or \textit{has no gap at $z\in C$ along any ray starting from $z$}, 
if for every $w\in C$ one has
$$
f(z)\le \limsup_{t\searrow 0} f(z+t(w-z)).
$$
Notice that this property is weaker than radial lower semicontinuity of $f$ at
$z\in C$.
\medbreak
A point $z\in C$ is called a \textit{flat point of $C$} if there are
three points $u,v, w$ in $C$ such that $z\in {]}u,v{[}$ and
there is no $w'\in C$ such that $z\in{]}w,w'{[}$.
The first condition means that $z$ is not an extreme point of $C$,
the second condition means that $z$ lies on the relative boundary of $C$
(or that $z$ is not in the intrinsic core of $C$).
Simple examples of convex sets with flat points are \textit{$n$-simplices}
$\Delta_n$,
that is convex hulls of $n+1$ affinely independent points in $E$:
any point on the boundary of $\Delta_n$ except the vertices is a flat point.
A convex set with no flat points is sometimes called a \textit{strictly convex set}.
Simple examples of strictly convex subsets of $E$ are affine subspaces,
finitely open convex subsets (i.e.\ convex sets $C$ such that
$C\cap F$ is open in $F$ for every finite dimensional subspace $F$ of $E$),
line segments or closed balls with respect to a
strictly convex norm. 

\begin{theorem}
Let $E$ be a real vector space with dual space $E^*$,
and let $C\subset E$ be a convex subset with more than one point.
Let $f:C\rightarrow\R$ be radially lower stable at
each flat point of $C$.
Then $f$ is convex if and only if there exists $c^*\in E^*$ not constant on $C$
such that for every $\lambda\in\R$, the function $f+\lambda c^*$ is
quasiconvex.
\end{theorem}
\medskip\noindent
\textbf{Proof.}
The necessity follows from the fact that $C$ has at least two points, so
there exists a linear form $c^*\in E^*$ not constant on $C$
by the basis extension theorem,
and the sum of two convex functions is convex, so
the function $f+\lambda c^*$ is quasiconvex.
We now prove the sufficiency.
Let $u, v\in C$, $u\ne v$, and let $t\in {]}0, 1{[}$. Set $z_t:=v+t(u-v)$.
We must show that
$$
f(z_t)\leq f(v)+t(f(u)-f(v)).
$$
First assume $\langle c^*, u-v\rangle\ne 0$.
Then we can find $\lambda\in\R$ such that
$$
\lambda\langle c^*, u-v\rangle=f(v)-f(u).
$$
Since $f+\lambda c^*$ is quasiconvex and
$(f+\lambda c^*)(u)=(f+\lambda c^*)(v)$, we derive that
\begin{align*}
(f+ \lambda c^*)(v) \ge (f+ \lambda c^*)(z_t) &=f(z_t)+\la \lambda c^*,v+t (u-v)\ra\\
&=f(z_t)+\la \lambda c^*,v\ra+t(f(v)-f(u)),
\end{align*}
which gives 
$$
f(z_t)\leq f(v)+t(f(u)-f(v)),
$$
as required.
\smallbreak
Now assume $\langle c^*, u-v\rangle= 0$.
Since $c^*$ is not constant on $C$, we may choose $w\in C$ such that
$\la c^*, w\ra\ne\la c^*, u\ra=\la c^*, v\ra$.
For $s\in {]}0, 1{[}$, consider the point
$v_s:=v+s(w-v)$
on the segment $[w,v]$ and the point $z_{t_s}:=v_s+t_s(u-v_s)$ at the
intersection of the segments $[u,v_s]$ and $[w,z_t]$; see the picture.
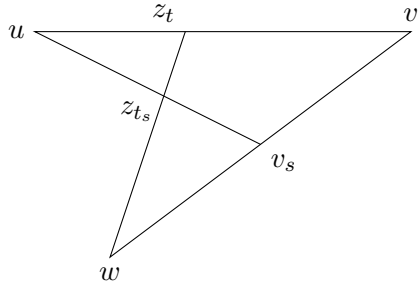
\pict
Since $\la c^*, w\ra\ne \la c^*, v\ra$,
we may apply the first part of the proof
with the points $w, v\in C$ and $v_s=v+s(w-v)$ to get
$$
f(v_s)\leq f(v)+s (f(w)-f(v)),
$$
and
since $\la c^*, v_s\ra= \la c^*, v\ra+s \la c^*, w-v\ra\ne \la c^*, v\ra=\la c^*, u\ra$,
we may also apply the first part of the proof
with the points $u, v_s\in C$ and $z_{t_s}=v_s+t_s(u-v_s)$ to get
\begin{align*}
f(z_{t_s})\leq (1-t_s)f(v_s)+t_s f(u).
\end{align*}
Combining these two inequalities, we derive that
\begin{align}\label{ineq}
f(z_{t_s})\leq (1-t_s)[f(v)+s (f(w)-f(v))]+t_s f(u).
\end{align}
As $s\rightarrow 0^+$, we have $t_s\to t$ and
the right-hand side of \eqref{ineq} tends to $(1-t)f(v)+t f(u)$.
On the other hand, $z_{t_s}\to z_t$ on the segment $[w,z_t]$.
Two cases are possible for the point $z_t$.
If $z_t$ is a flat point of $C$, then by the stability assumption
$$
f(z_t)\le \limsup_{s\searrow 0} f(z_{t_s}).
$$
If $z_t$ is not a flat point of $C$, we can find a point
$w'\in C$ such that $z_t\in {]}w,w'{[}$. Clearly,
$\la c^*, w'-w\ra\ne 0$, because $w'-w$ is a non-zero
multiple of $z_t-w$ and $\la c^*, z_t-w\ra\ne 0$ since
$\la c^*, z_t\ra=\la c^*, u\ra\ne \la c^*, w\ra$.
Thus, by the first part of the proof, $f$ is convex on $[w,w']$.
But a finite convex function on an interval is automatically continuous
on the interior of that interval. Since  $z_t\in {]}w,w'{[}$,
it follows that
$$
f(z_t)= \lim_{s\searrow 0} f(z_{t_s}).
$$
Therefore, in both cases, letting $s\rightarrow 0^+$ in \eqref{ineq},
we obtain
$$
f(z_t)\leq f(v)+t(f(u)-f(v)).
$$
The proof is complete. $\hfill\Box$
\medbreak\noindent
{\it Remarks.}
1. If $C$ has no flat points (that is, $C$ is strictly convex),
the regularity assumption on $f$ is automatically satisfied.
Otherwise, this assumption cannot be dropped.
Indeed, in $E=\R^2$,
consider the triangle $C:={\rm conv\,}\{u,v,w\}$ with vertices
$u=(1,0)$, $v=(0,1)$ and $w=(0,0)$,
and define $c^*\in E^*$ by
$$c^*:x:=(x_1,x_2)\in \R^2\mapsto \la c^*,x\ra:=x_1+x_2.$$
Then, for every $x\in C$ and $z\in [u,v]$ one has
$
0\le \la c^*,x\ra\le \la c^*,z\ra=1.
$

Now let $f:C\to\R$ given by $f (x) = 1$
if $x\in [u,v{[}$, $f (x) = 0$ otherwise.
This function is not convex since it is not convex on $[u,v]$.
But $g:=f +\lambda c^*$ is quasiconvex for all $\lambda\in\R$.
Indeed, $g$ being equal to $1+\lambda c^*$ on $[u,v{[}$ and to $\lambda c^*$
on $C\setminus [u,v{[}$, is (quasi)convex on each of these convex subsets.
Now, let  $y\in {]}x,z{[}$ with $x\in C\setminus [u,v{[}$ and $z\in [u,v{[}$.
We show that
$g(y)\le \max \{g(x),g(z)\}.$ 
If $y\in [u,v{[}$ (which happens if $x=v$), then
$g(y)=g(z)\le \max \{g(x),g(z)\}$. Otherwise, $g(y)= \lambda \la c^*,y\ra$.
Since $\la c^*,x\ra\le \la c^*,z\ra$ and $y\in {]}x,z{[}$,
it follows that $\la c^*,x\ra\le \la c^*,y\ra\le \la c^*,z\ra$.
Hence, for any $\lambda\in \R$,
\begin{align*}
g(y)=\lambda\langle c^*, y\rangle &\leq \max\{\lambda\langle c^*, x\rangle,\lambda\langle c^*, z\rangle\}\\
&\leq\max\{\lambda\langle c^*, x\rangle,1+\lambda\langle c^*, z\rangle\}\\
&= \max\{g(x), g(z)\}.
\end{align*}
This shows that $g=f +\lambda c^*$ is quasiconvex for all $\lambda\in\R$.
Incidentally, one easily check that $f$ is not radially lower
stable at any of the flat points $z\in {]}u,v{[}$ because one has
$\limsup_{t\searrow 0} f(z+t(w-z))=0< f(z)=1$.
\smallbreak

2. The assumption that \textit{$c^*$ is not constant on $C$} cannot be
omitted since for any quasiconvex function $f$ and any $\lambda\in \R$
the function $u\mapsto f(u) +\lambda$ is quasiconvex.
\smallbreak

3. The assumption
\begin{align*}
f+\lambda c^* \textit{ is quasiconvex for every } \lambda\in\R
\end{align*}
cannot be relaxed to
\begin{align*}
f+\lambda c^* \textit{ is quasiconvex for every } \lambda\ge 0
\ (\textit{or every } \lambda\le 0).
\end{align*}
Indeed, in $E=\R^2$,
consider the same triangle $C={\rm conv\,}\{u,v,w\}$ and linear form
$c^*\in E^*$ as in Remark 1 above, and define $f:C\to\R$ by
$f (x) = 0$ if $x\in {]}u,v{[}$, $f (x) = 1$ otherwise.
This $f$ is radially lower stable at
any flat point of $C$ and $f+\lambda c^*$ is quasiconvex for every 
$\lambda\le 0$, yet $f$ is not convex.

Similarly, the above assumption cannot be relaxed to
\begin{align*}
f+\lambda c^* \textit{ is quasiconvex for every {\em arbitrarily small} }
\lambda\in\R.
\end{align*}
Indeed, there exist non-convex quasiconvex functions satisfying
such a property, see \cite{BGJ12,PhuAn96}.
\smallbreak

4.
An analogue property linking monotone and quasi-monotone operators was
established by Hadjisavvas \cite{Had06}.
For the relationships between (quasi)convex functions and (quasi)monotone operators,
see for instance \cite{ACL94}.
\medbreak\noindent
\textit{Acknowledgement.}
This work was completed while the authors were visiting the
Vietnam Institute for Advanced Study in Mathematics (VIASM).
They would like to thank the VIASM for financial support and hospitality.
The first author is funded by
Vietnam National Foundation for Science and Technology Development (NAFOSTED)
under grant number 101.01-2014.56.
The authors also gratefully acknowledge the anonymous referees for
their helpful remarks
that allowed a substantial improvement of the presentation.


{\small
\begin{thebibliography}{9}
\bibitem{ACL94} {D. Aussel, J.-N. Corvellec and M. Lassonde}:
{\em Subdifferential characterization of quasiconvexity and convexity}, J. Convex Anal. \textbf{1} (1994), 195--201.

\bibitem{BGJ12} {E.N. Barron, R. Goebel and R.R. Jensen}:
{\em Functions which are quasiconvex under linear perturbations},
SIAM J. Optim. \textbf{22} (2012), 1089--1108.

\bibitem{Cro77} {J.-P. Crouzeix}: {\rm Contributions \`a l'\'etude des
fonctions quasi-convexes}, Th\`ese d'\'Etat, Universit\'e de Clermont-Ferrand II, 1977.

\bibitem{Had06} {N. Hadjisavvas}: {\em Translations of quasimonotone maps and monotonicity},
Appl. Math. Lett. \textbf{19} (2006), 913--915.

\bibitem{PhuAn96} {H.X. Phu and P.T. An}:
{\em Stable generalization of convex functions},
Optimization {\bf 38} (1996), 309--318.

\end{thebibliography}
}
\end{document}